\documentclass[a4paper,12pt]{amsart}
\usepackage{amssymb}
\usepackage{ifthen}
\usepackage{graphicx}
\usepackage{float}
\newtheorem{thm}{Theorem}
\newtheorem{cor}{Corollary}
\newtheorem{lem}{Lemma}

\newtheorem{que}{Question}
\newtheorem{rem}{Remark}
\newtheorem{prob}{Problem}
\newtheorem{conj}{Conjecture}
\theoremstyle{definition}
\newtheorem{example}[equation]{Example}%[section]

\newcounter {own}
\def\theown {\thesection       .\arabic{own}}

\newenvironment{pf}[1][]{%
 \vskip 3mm
 \noindent
 \ifthenelse{\equal{#1}{}}%
  {{\slshape Proof. }}%
  {{\slshape #1.} }%
 }%
{\qed\bigskip}

\newcounter{alphabet}
\newcounter{tmp}

\makeatletter
\newcommand{\Ref}[1]{\@ifundefined{r@#1}{}{\setcounter{tmp}{\ref{#1}}\Alph{tmp}}}
\makeatother
\newcounter{minutes}
\divide\time by 60
\newcounter{hours}
\multiply\time by 60
\addtocounter{minutes}{-\time}
\newcommand{\IN}{{\mathbb N}}
\newcommand{\IC}{{\mathbb C}}
\newcommand{\ID}{{\mathbb D}}

\def\be{\begin{equation}}
\def\ee{\end{equation}}

\newcommand{\bee}{\begin{enumerate}}
\newcommand{\eee}{\end{enumerate}}

\newcommand{\blem}{\begin{lem}}
\newcommand{\elem}{\end{lem}}
\newcommand{\bthm}{\begin{thm}}
\newcommand{\ethm}{\end{thm}}
\newcommand{\bcor}{\begin{cor}}
\newcommand{\ecor}{\end{cor}}
\newcommand{\beg}{\begin{example}}
\newcommand{\eeg}{\end{example}}
\newcommand{\bques}{\begin{que}}
\newcommand{\eques}{\end{que}}
\newcommand{\begs}{\begin{examples}}
\newcommand{\eegs}{\end{examples}}
\newcommand{\bdefe}{\begin{defin}}
\newcommand{\edefe}{\end{defin}}
\newcommand{\bprob}{\begin{prob}}
\newcommand{\eprob}{\end{prob}}
\newcommand{\bei}{\begin{itemize}}
\newcommand{\eei}{\end{itemize}}

\newcommand{\bcon}{\begin{conj}}
\newcommand{\econ}{\end{conj}}
\newcommand{\bcons}{\begin{conjs}}
\newcommand{\econs}{\end{conjs}}
\newcommand{\bprop}{\begin{propo}}
\newcommand{\eprop}{\end{propo}}
\newcommand{\br}{\begin{rem}}
\newcommand{\er}{\end{rem}}
\newcommand{\brs}{\begin{rems}}
\newcommand{\ers}{\end{rems}}
\newcommand{\bo}{\begin{obser}}
\newcommand{\eo}{\end{obser}}
\newcommand{\bos}{\begin{obsers}}
\newcommand{\eos}{\end{obsers}}
\newcommand{\bpf}{\begin{pf}}
\newcommand{\epf}{\end{pf}}
\newcommand{\ba}{\begin{array}}
\newcommand{\ea}{\end{array}}
\newcommand{\beq}{\begin{eqnarray}}
\newcommand{\beqq}{\begin{eqnarray*}}
\newcommand{\eeq}{\end{eqnarray}}
\newcommand{\eeqq}{\end{eqnarray*}}

\newcommand{\ds}{\displaystyle}

\def\cc{\setcounter{equation}{0}   % THIS CLEARS THE COUNTER
\setcounter{figure}{0}\setcounter{table}{0}}

\begin{document}

\bibliographystyle{amsplain}

%
%\begin{center}
%{\tiny \texttt{FILE:~\jobname .tex,
%        printed: \number\year-\number\month-\number\day,
%        \thehours.\ifnum\theminutes<10{0}\fi\theminutes}
%}
%\end{center}

\title[Logarithmic Coefficients and a Coefficient Conjecture for Univalent
Functions]{Logarithmic Coefficients and a Coefficient Conjecture for Univalent Functions}
%{A lemma of rogosinski and Mertens function}

%=========================================================================
\thanks{%$^\dagger$
File:~\jobname .tex,
          printed: \number\day-\number\month-\number\year,
          \thehours.\ifnum\theminutes<10{0}\fi\theminutes}
%=========================================================================

\author[M. Obradovi\'{c}]{Milutin Obradovi\'{c}}
\address{M. Obradovi\'{c},
Department of Mathematics,
Faculty of Civil Engineering, University of Belgrade,
Bulevar Kralja Aleksandra 73, 11000
Belgrade, Serbia. }
\email{obrad@grf.bg.ac.rs}

\author[S. Ponnusamy]{Saminathan Ponnusamy
%$^\dagger $
%${}^{~\mathbf{*}}$
}
\address{S. Ponnusamy, Indian Statistical Institute (ISI), Chennai Centre, SETS (Society
for Electronic Transactions and Security), MGR Knowledge City, CIT
Campus, Taramani, Chennai 600 113, India.}
\email{samy@isichennai.res.in, samy@iitm.ac.in}

\author[K.-J. Wirths]{Karl-Joachim Wirths}
\address{K.-J. Wirths, Institut f\"ur Analysis und Algebra, TU Braunschweig,
38106 Braunschweig, Germany.}
\email{kjwirths@tu-bs.de}

\subjclass[2010]{30C45}
\keywords{Univalent, starlike, convex and close-to-convex functions, subordination, logarithmic  coefficients
and coefficient estimates
%\\
%$%{}^{\mathbf{*}}
%^\dagger$ {\tt This author is on leave from the Department of Mathematics,
%Indian Institute of Technology Madras, Chennai-600 036, India}
}
%\date{\today  %June. 30, 09
%;  File: 2013(2).tex}

\begin{abstract}
Let ${\mathcal U}(\lambda)$ denote the family of analytic functions $f(z)$, $f(0)=0=f'(0)-1$, in the unit disk $\ID$, which satisfy the condition
$\big |\big (z/f(z)\big )^{2}f'(z)-1\big |<\lambda $ for some $0<\lambda \leq 1$. The logarithmic coefficients $\gamma_n$ of $f$ are defined by the formula
$\log(f(z)/z)=2\sum_{n=1}^\infty \gamma_nz^n$. In a recent paper, the present authors proposed a conjecture that if $f\in {\mathcal U}(\lambda)$ for
some $0<\lambda \leq 1$,  then $|a_n|\leq \sum_{k=0}^{n-1}\lambda ^k$ for $n\geq 2$ and provided a new proof for the case $n=2$.
One of the aims of this article is to present a proof of this conjecture for $n=3, 4$ and an elegant proof of the inequality for $n=2$, with
equality for $f(z)=z/[(1+z)(1+\lambda z)]$.
In addition, the authors prove the following sharp inequality for $f\in{\mathcal U}(\lambda)$:
$$\sum_{n=1}^{\infty}|\gamma_{n}|^{2} \leq \frac{1}{4}\left(\frac{\pi^{2}}{6}+2{\rm Li\,}_{2}(\lambda)+{\rm Li\,}_{2}(\lambda^{2})\right),
$$
where ${\rm Li}_2$ denotes the dilogarithm function.
Furthermore, the authors prove two such new inequalities satisfied by the corresponding
logarithmic coefficients of some other subfamilies of $\mathcal S$.
\end{abstract}

%\thanks{The work of the second author was supported by MNZZS Grant, No. ON174017, Serbia. }

\maketitle
\pagestyle{myheadings}
\markboth{M. Obradovi\'{c}, S. Ponnusamy and K.-J. Wirths}{Logarithmic  Coefficients for Univalent Functions}
\cc

\section{Introduction}
Let $\mathcal A$ be the class of functions $f$  analytic in the unit
disk $\ID=\{z\in \IC:\, |z|<1\}$ with the normalization $f(0)=0=f'(0)-1$. Let $\mathcal S$ denote the class of functions
$f$ from $\mathcal A$ that are univalent in $\ID$. Then the logarithmic coefficients $\gamma_n$ of $f\in {\mathcal S}$ are defined by the formula
\begin{equation}\label{eq:log}
\frac{1}{2}\log\left (\frac{f(z)}{z}\right )=\sum_{n=1}^\infty \gamma_nz^n, \quad z\in \ID.
\end{equation}
These coefficients play an important role for various estimates in the theory of univalent functions.
When we require a distinction,  we use the notation $\gamma_n(f)$  instead of $\gamma_n$. For example, the Koebe
function $k(z)=z(1-e^{i\theta}z)^{-2}$ for each $\theta$ has logarithmic coefficients $\gamma _n(k)=e^{in\theta}/n$, $n\geq 1$. If
$f\in {\mathcal S}$ and $f(z)=z+\sum_{n=2}^{\infty} a_{n}z^{n}$, then by \eqref{eq:log} it follows that $2\gamma_1=a_2$ and hence,
by the Bieberbach inequality, $|\gamma _1|\leq 1$. Let ${\mathcal S}^{\star}$ denote the class of functions $f\in{\mathcal S}$ such that $f(\ID)$ is
starlike with respect to the origin. Functions $f\in{\mathcal S}^{\star}$ are characterized by the condition ${\rm Re}\, (zf'(z)/f(z))>0$ in $\ID$.
The inequality $|\gamma _n|\leq 1/n$ holds for starlike functions $f\in {\mathcal S}$, but is false for the
full class $\mathcal S$, even in order of magnitude. See \cite[Theorem 8.4 on page~242]{Duren}.  In \cite{Girela00}, Girela pointed out that
this bound is actually false for the class of close-to-convex functions in $\ID$ which is defined as follows:
A function $f\in\mathcal{A}$ is called close-to-convex, denoted by $f\in\mathcal{K}$, if there exists a real $\alpha$ and a $g\in  {\mathcal S}^{\star}$
such that
$${\rm Re} \left( e^{i\alpha}\frac{zf'(z)}{g(z)} \right) > 0, \quad \mbox{$z\in \ID$.}
$$
For $0\leq \beta <1$, a function $f\in {\mathcal S}$ is said to belong to the class of starlike functions of order $\beta$,
denoted by $f\in {\mathcal S}^{\star}(\beta)$, if ${\rm Re} \left (zf'(z)/f(z)\right )>\beta $ for $z\in\ID$.
Note that ${\mathcal S}(0)=:{\mathcal S}^{\star}$. The class of all convex functions of order $\beta$, denoted by ${\mathcal C}(\beta)$,
is then defined by ${\mathcal C}(\beta)=\{f\in\mathcal{S}:\, zf'\in {\mathcal S}^{\star}(\beta)\}$.
The class ${\mathcal C}(0)=:{\mathcal C}$ is usually referred to as the class of convex functions in $\ID$.
With the class $\mathcal S$ being of the first priority, its subclasses such as ${\mathcal S}^{\star}$, $ {\mathcal K}$,
and $\mathcal C$, respectively, have been extensively studied in the literature and they appear in
different contexts.
We refer to  \cite{Duren,Go,OPW,P} for a general reference related to the present study.
In \cite[Theorem 4]{DL79}, it was shown that the logarithmic coefficients $\gamma_n$ of every function $f\in {\mathcal S}$ satisfy
\be\label{eq36}
\sum_{n=1}^{\infty}|\gamma_{n}|^{2}\leq\frac{\pi^{2}}{6}
\ee
and the equality is attained for the  Koebe function.
The proof uses ideas from the work of Baernstein \cite{Baerns74} on integral means. However, this result is easy to prove (see Theorem \ref{13-th 12})
in the case of functions in the class $\mathcal{U}:=\mathcal{U}(1)$ which is defined as follows:
$$ \mathcal{U}(\lambda)=\left \{f\in {\mathcal A}:\, \left|\left(\frac{z}{f(z)}\right)^2f'(z)-1\right|< \lambda, \quad z\in \ID \right \},
$$
where $\lambda \in (0,1]$. It is known that \cite{Aks58,AksAvh70,OzNu72} every $f\in \mathcal{U}$ is  univalent in $\ID$ and hence,
$\mathcal{U}(\lambda)\subset \mathcal{U} \subset \mathcal{S} $ for $\lambda \in (0,1]$.  The present authors have
established many interesting properties of the family $\mathcal{U}(\lambda)$. See \cite{OPW} and the references therein. For example, if $f\in {\mathcal U}(\lambda)$ for some $0<\lambda \leq 1$ and
$a_2 =f''(0)/2$, then we have the subordination relations
\be\label{OPW7-eq4c}
\frac{f(z)}{z}\prec \frac{1}{1+(1+\lambda)z+\lambda z^2}=\frac{1}{(1+z)(1+\lambda z)}, ~ z\in \ID,
\ee
and
$$ \frac{z}{f(z)}+ a_2z \prec 1+2\lambda z + \lambda z^2, ~ z\in \ID.
$$
Here $\prec$ denotes the usual subordination  \cite{Duren,Go,P}. In addition, the following conjecture was proposed in \cite{OPW}.

\bcon\label{conj1}
Suppose that $f\in {\mathcal U}(\lambda)$ for some $0<\lambda \leq 1$.  Then $|a_n|\leq \sum_{k=0}^{n-1}\lambda ^k$ for $n\geq 2$.
\econ

In Theorem \ref{13-th 12b}, we present a direct proof of an inequality analogous to \eqref{eq36} for functions in ${\mathcal U}(\lambda)$
and in Corollary \ref{13-th 12b}, we obtain the inequality \eqref{eq36} as a special case for ${\mathcal U}$. At the end of Section \ref{sec2},
we also consider estimates of the type \eqref{eq36}
for some interesting subclasses of univalent functions. However, Conjecture \ref{conj1} remains open for $n\geq 5$. On the other hand, the proof
for the case $n=2$ of this conjecture is due to \cite{VY2013} and an alternate proof was
obtained recently by the present authors in \cite[Theorem 1]{OPW}. In this paper, we show that Conjecture \ref{conj1} is true for $n=3, 4$.
and our proof includes an elegant proof of the case $n=2$. The main results and their proofs are presented in Sections \ref{sec2} and \ref{sec3}.

\section{Logarithmic coefficients of functions in ${\mathcal U}(\lambda)$\label{sec2}}

\bthm\label{13-th 12b}
For $0<\lambda \leq 1$, the logarithmic coefficients of $f\in\mathcal{U}(\lambda)$ satisfy the inequality
\be\label{eq36b}
\sum_{n=1}^{\infty}|\gamma_{n}|^{2} \leq \frac{1}{4}\left(\frac{\pi^{2}}{6}+2{\rm Li\,}_{2}(\lambda)+{\rm Li\,}_{2}(\lambda^{2})\right),
\ee
where ${\rm Li}_2$ denotes the dilogarithm function given by
$${\rm Li}_2(z)= \sum_{n=1}^\infty \frac{z^{n}}{n^2}=  z\int_0^1\frac{\log (1/t)}{1-tz}\,dt.
$$
The inequality \eqref{eq36b} is sharp. Further, there exists a function $f \in\mathcal{U}$ such that $|\gamma_{n}|>(1+\lambda ^n)/(2n)$ for some $n$.
\ethm
\bpf Let  $f\in \mathcal{U}(\lambda)$. Then, by \eqref{OPW7-eq4c},  we have
$$\frac{z}{f(z)}\prec (1-z)(1-\lambda z)
$$
which clearly gives
\be\label{eq36c}
\sum_{n=1}^{\infty} \gamma_{n}z^{n}= \log\sqrt{\frac{f(z)}{z}}\prec\frac{ -\log (1-z)-\log (1-\lambda z)}{2}=
\sum_{n=1}^{\infty} \frac{1}{2n}(1+\lambda ^{n})z^{n}.
\ee
Again, by Rogosinski's theorem (see \cite[6.2]{Duren}), we obtain
$$\sum_{n=1}^{\infty}|\gamma_{n}|^{2} \leq \sum_{n=1}^{\infty}\frac{1}{4n^{2}}(1+\lambda ^{n})^{2}
=\frac{1}{4}\left(\sum_{n=1}^{\infty}\frac{1}{n^{2}}+2\sum_{n=1}^{\infty}\frac{\lambda ^{n}}{n^{2}}+
\sum_{n=1}^{\infty}\frac{\lambda ^{2n}}{n^{2}}\right)
$$
and the desired inequality \eqref{eq36b} follows. For the function
$$g_{\lambda}(z)=\frac{z}{(1-z)(1-\lambda z)},
$$
we find that $\gamma_{n}(g_{\lambda}) =(1+\lambda ^n)/(2n)$ for $n\ge 1$ and therefore, we have the equality in \eqref{eq36b}.
Note that $g_{1}(z)$ is the Koebe function $z/(1-z)^2$.

\begin{figure}%[H]
\begin{center}
\includegraphics[height=5.5cm, width=5.5cm, scale=1]{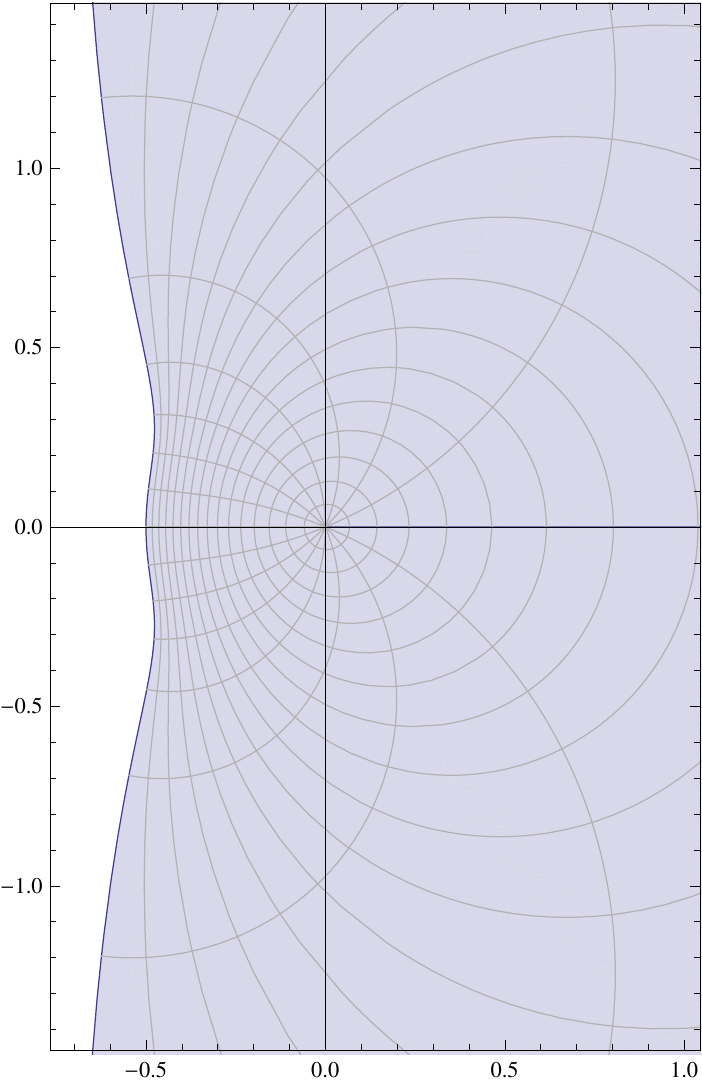}
\hspace{1cm}
\includegraphics[height=5.5cm, width=5.5cm, scale=1]{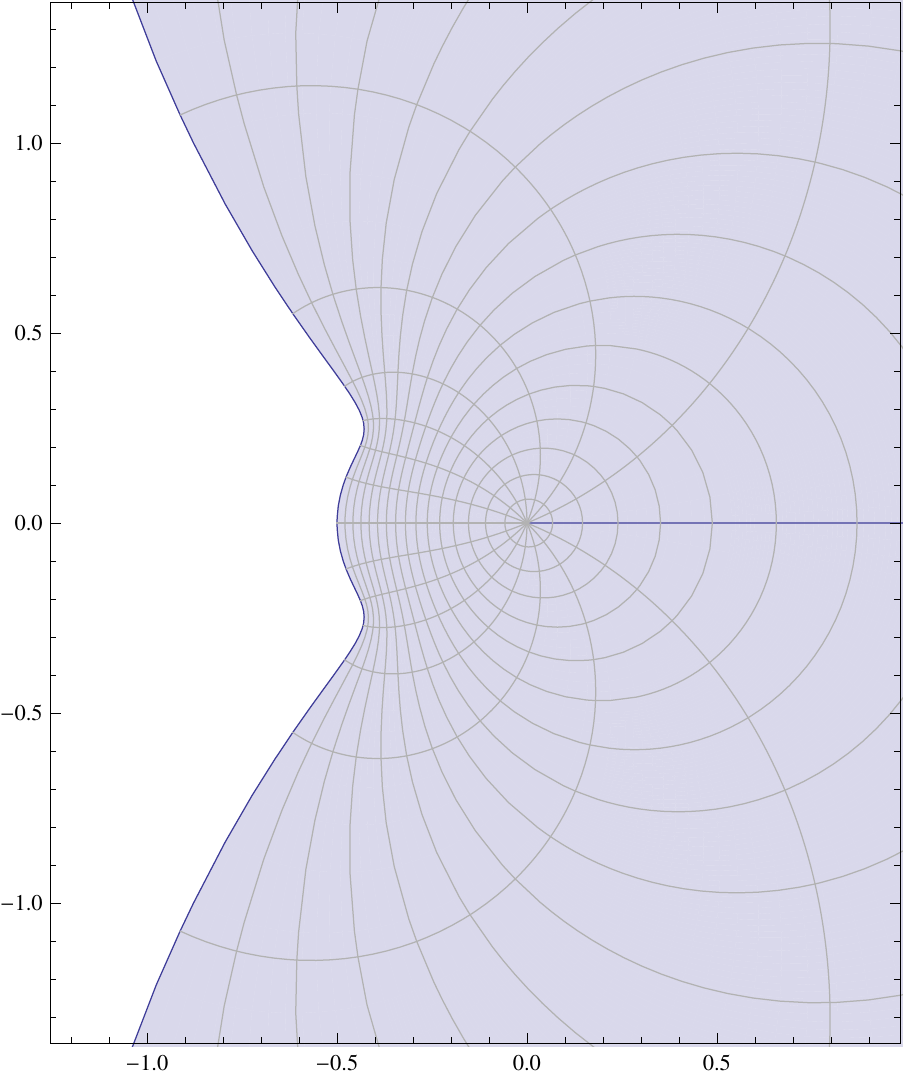}
\end{center}
(a )$\lambda=0.25$ \hspace{6cm} (b) $\lambda=0.5$

\vspace{0.25cm}

\begin{center}
\includegraphics[height=5.5cm, width=5.5cm, scale=1]{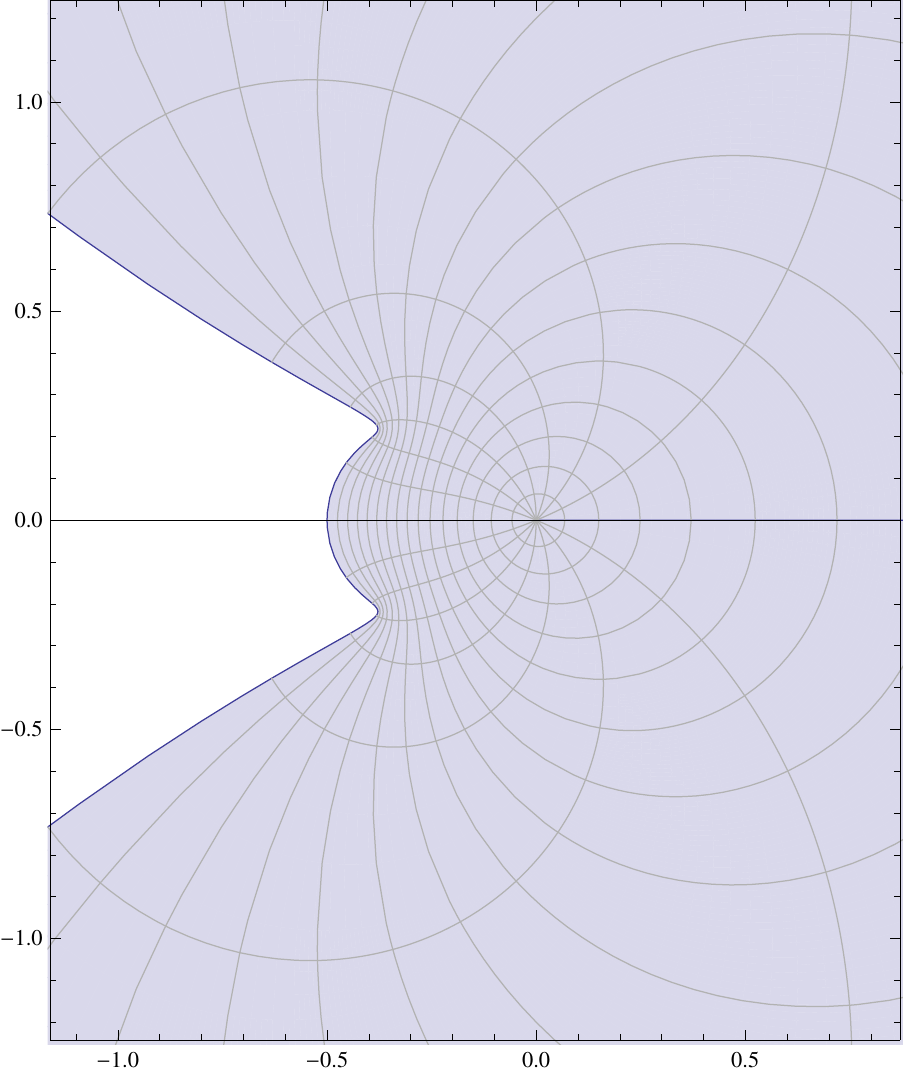}
\hspace{1cm}
\includegraphics[height=5.5cm, width=5.5cm, scale=1]{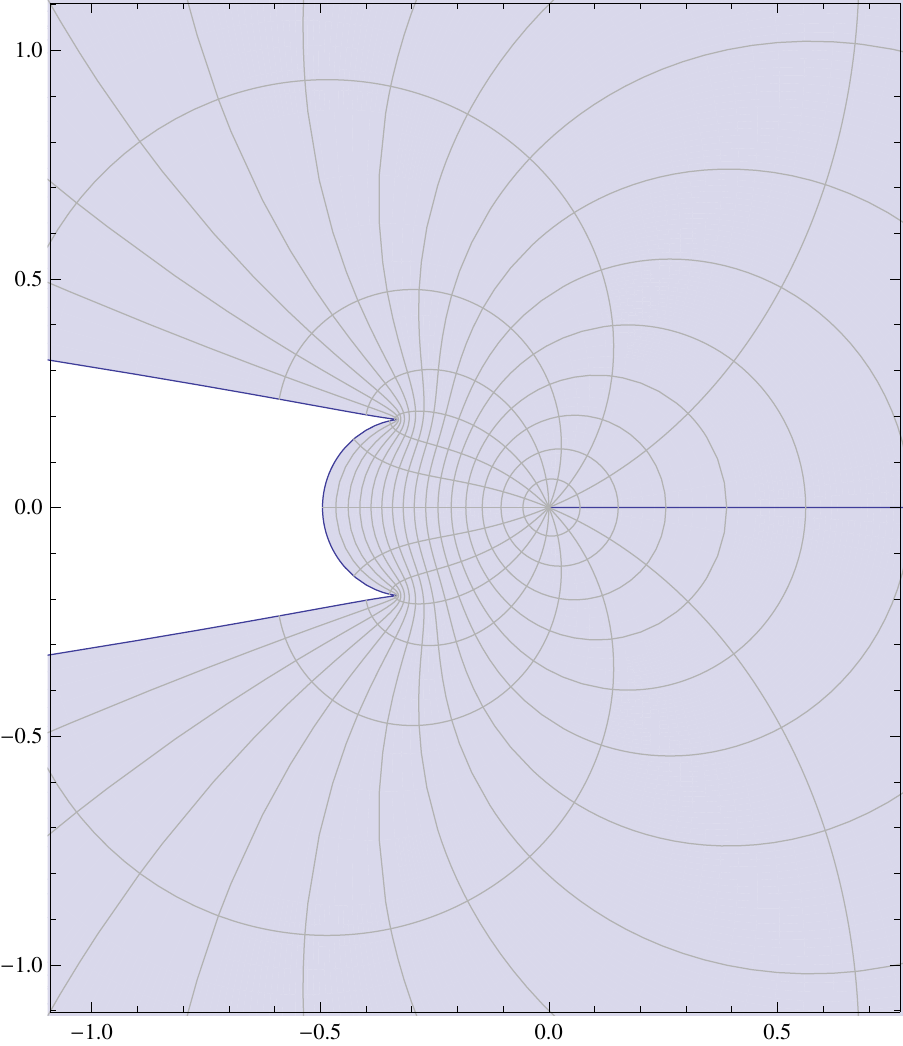}
\end{center}
(c) $\lambda=0.75$ \hspace{6cm} (d) $\lambda=1$

\caption{The image of $\ds f_{\lambda}(z)= \frac{z}{(1-z)(1-\lambda z)(1+(\lambda/(1+\lambda)) z)}$  under $\ID$ for certain values of $\lambda$
\label{logCoeffig1}}
\end{figure}

From the relation \eqref{eq36c},  we cannot conclude that
$$|\gamma_{n}(f)|\leq |\gamma_{n}(g_{\lambda})|=\frac{1+\lambda ^n}{2n} ~\mbox{ for $f\in\mathcal{U}(\lambda)$}.
$$
Indeed for the function $f_{\lambda}$ defined by
\be\label{eqext1}
f_{\lambda}(z)= \frac{z}{(1-z)(1-\lambda z)(1+(\lambda/(1+\lambda)) z)}
\ee
we find that
$$\frac{z}{f_{\lambda}(z)}= 1+ \frac{\lambda -(1+\lambda)^2}{1+\lambda} z +\frac{\lambda ^2}{1+\lambda}z^3
$$
and
$$\left(\frac{z}{f_{\lambda}(z)}\right)^2f_{\lambda}'(z)-1= -\frac{2\lambda ^2}{1+\lambda}z^3 =-\left (1- \frac{(1+2\lambda)(1-\lambda )}{1+\lambda}\right )z^3
$$
which clearly shows that $f_{\lambda}\in \mathcal{U}(\lambda)$.
The images of $f_{\lambda}(z)$ under $\ID$ for certain values of $\lambda$ are shown in Figures \ref{logCoeffig1}(a)-(d).
Moreover, for this function,  we have
\beqq
\log \left (\frac{f_{\lambda}(z)}{z}\right )&=&-\log(1-z)-\log (1-\lambda z)-\log\left(1+\frac{\lambda}{1+\lambda} z\right)\\
&=& 2\sum_{n=1}^{\infty} \gamma_{n}(f_{\lambda})z^n,
%\quad \gamma_{n}(f_{\lambda}) =\frac{1}{2}\left (\frac{1+\lambda ^{n}}{n} + (-1)^{n}\frac{\lambda ^n}{(1+\lambda )^n}\right ).
\eeqq
where
$$\gamma_{n}(f_{\lambda}) =\frac{1}{2}\left (\frac{1+\lambda ^{n}}{n} + (-1)^{n}\frac{\lambda ^n}{n(1+\lambda )^n}\right ).
$$
This contradicts the above inequality at least for even integer values of $n\geq 2$. Moreover, with these $\gamma_{n}(f_{\lambda})$
for $n\geq 1$,  we obtain
\beqq
\sum_{n=1}^{\infty}|\gamma_{n}(f_{\lambda})|^{2} &= &  \frac{1}{4}  \sum_{n=1}^{\infty}
 \left \{\frac{(1+\lambda ^{n})^2}{n^2}+  2\frac{(-1)^{n}}{n^2} \left [ \left (\frac{\lambda ^{2} }{1+\lambda}\right ) ^{n} \right. \right. \\
& & \left .\left . + \left (\frac{\lambda}{1+\lambda}\right )^{n} \right ]  +\frac{1}{n^2} \left (\frac{\lambda}{1+\lambda}\right )^{2n} \right \}
%$& =  \frac{1}{4} (I_1+I_2+I_3). %=\frac{\pi^{2}}{6}+\frac{1}{12}-\log\frac{3}{2}<\frac{\pi^{2}}{6}
\eeqq
and by a computation, it follows easily that
\beqq
\sum_{n=1}^{\infty}|\gamma_{n}(f_{\lambda})|^{2}& =&\frac{1}{4}\left(\frac{\pi^{2}}{6}+2{\rm Li\,}_{2}(\lambda)+{\rm Li\,}_{2}(\lambda^{2})\right) \\
&& \qquad +\frac{1}{2} \left[ {\rm Li\,}_{2}\left (\frac{-\lambda ^{2} }{1+\lambda}\right )
+ {\rm Li\,}_{2}\left (\frac{-\lambda }{1+\lambda}\right )\right]
  +\frac{1}{4}{\rm Li\,}_{2}\left (\frac{\lambda^{2}}{(1+\lambda)^{2}}\right ) \\
  & =&\frac{1}{4}\left(\frac{\pi^{2}}{6}+2{\rm Li\,}_{2}(\lambda)+{\rm Li\,}_{2}(\lambda^{2})\right)  +\frac{1}{4} A(\lambda ) \\
&<& \frac{1}{4}\left(\frac{\pi^{2}}{6}+2{\rm Li\,}_{2}(\lambda)+{\rm Li\,}_{2}(\lambda^{2})\right) ~\mbox{ for $0<\lambda \leq 1$,}
\eeqq
and we complete the proof, provided $A(\lambda )<0$ for $0<\lambda \leq 1$.
%Here
%\be\label{MOFM-eq1}
%{\rm Li\,}_{2}(z) =\sum_{n=1}^{\infty} \frac{z^n}{n^2} =\int_{0}^{1}\frac{z}{1-tz}\log (1/t)\, dt.
%\ee
Now, we claim that
$$A(\lambda ):=2\left[ {\rm Li\,}_{2}\left (\frac{-\lambda ^{2} }{1+\lambda}\right )
+ {\rm Li\,}_{2}\left (\frac{-\lambda }{1+\lambda}\right )\right]
  +{\rm Li\,}_{2}\left (\frac{\lambda^{2}}{(1+\lambda)^{2}}\right ) <0.
$$
Because ${\rm Li\,}_{2}(z^2)=2({\rm Li\,}_{2}(z)+{\rm Li\,}_{2}(-z))$, the last claim is equivalent to
$$ \frac{A(\lambda)}{2}=2\,{\rm Li\,}_{2}\left (\frac{-\lambda }{1+\lambda}\right )
+ \left [ {\rm Li\,}_{2}\left (\frac{\lambda}{1+\lambda}\right )
+ {\rm Li\,}_{2}\left (\frac{-\lambda ^{2} }{1+\lambda}\right )\right] <0
$$
for $0<\lambda \leq 1$.  According to the integral representation of ${\rm Li\,}_{2}(z)$ given in the statement of
Theorem \ref{13-th 12b}, we can write
$$ A(\lambda) =-2\lambda \int_{0}^{1}B(\lambda, t) \log (1/t)\, dt,
$$
where
\beqq
B(\lambda, t)&=& \frac{2}{1+\lambda +t\lambda} -\frac{1}{1+\lambda -t\lambda} +\frac{\lambda }{1+\lambda +t\lambda ^2}\\
&=& \frac{(1+\lambda)- 3t\lambda }{(1+\lambda)^2 -t^2\lambda^2} +\frac{\lambda }{1+\lambda +t\lambda ^2}\\
&=& \frac{N(\lambda, t)}{[(1+\lambda)^2 -t^2\lambda^2][1+\lambda +t\lambda ^2]}
\eeqq
with
$$N(\lambda, t)=(1+\lambda)^3-(3-\lambda)(1+\lambda)\lambda t -4\lambda ^3t^2.
$$
Clearly, $B(1, t)>0$  for $t\in [0,1)$ and it follows that, $A(1)<0$. On the other hand, since $N(\lambda, t)$ is a decreasing function of
$t$ for $t\in [0,1]$, we obtain that
$$N(\lambda, t)\geq N(\lambda, 1)=(1+\lambda)^3-(3-\lambda)(1+\lambda) -4\lambda ^3= 1-\lambda^3+\lambda^2(1-\lambda)>0
$$
for $0<\lambda <1$. Consequently, $B(\lambda, t)>0$ for all $t\in [0,1]$ and for  $0<\lambda <1 $. This observation shows that $A(\lambda)<0$ for
$0<\lambda \leq 1 $. This proves the claim and thus, the proof is complete.
\epf

\bcor\label{13-th 12}
The logarithmic coefficients of $f\in\mathcal{U}$ satisfy the inequality
\be\label{eq36a}
\sum_{n=1}^{\infty}|\gamma_{n}|^{2}\leq \sum_{n=1}^{\infty}\frac{1}{n^{2}}=\frac{\pi^{2}}{6}.
\ee
We have equality in the last inequality for the Koebe function $k(z)=z(1-e^{i\theta}z)^{-2}$.
Further there exists a function $f \in\mathcal{U}$ such that $|\gamma_{n}|>1/n$ for some $n$.
\ecor

 \br
From the analytic characterization of starlike functions,  it is easy to see that for $f\in\mathcal{S}^{\star}$,
$$\frac{zf'(z)}{f(z)} -1=z\left (\log\left (\frac{f(z)}{z}\right )\right )' = 2\sum_{n=1}^{\infty} n \gamma_{n}z^{n} \prec \frac{2z}{1-z}
$$
and thus, by Rogosinski's result, we obtain that $|\gamma_{n}|\leq 1/n$ for $n\ge 1$. In fact for starlike functions of order $\alpha$, $\alpha \in [0,1)$,
the corresponding logarithmic coefficients  satisfy the inequality $|\gamma_{n}|\leq (1-\alpha)/n$ for $n\geq 1$. Moreover, one can quickly obtain that
$$\sum_{n=1}^{\infty}|\gamma_{n}|^{2}\leq (1-\alpha)^2 \frac{\pi^2}{6}
$$
if $f\in\mathcal{S}^{\star}(\alpha)$, $\alpha \in [0,1)$ (See also the proof of Theorem \ref{13-th 23b} and Remark \ref{rem3}).
As remarked in the proof of Theorem \ref{13-th 12b}, from the relation \eqref{eq36a}, we cannot conclude the same fact, namely,
$|\gamma_{n}|\leq 1/n$ for $n\ge 1$, for the class $\mathcal{U}$ although the Koebe function $k(z)=z/(1-z)^2$ belongs to
$\mathcal{U}\cap \mathcal{S}^{\star}$. For example,  if we set $\lambda =1$ in \eqref{eqext1}, then we have
$$\frac{z}{f_1(z)}=(1-z)^{2}\left(1+\frac{z}{2}\right)=1-\frac{3}{2}z+\frac{z^{3}}{2},
$$
where $f_1\in {\mathcal U}$ and for this function,  we obtain
\beqq
\sum_{n=1}^{\infty}|\gamma_{n}(f_1)|^{2}&= & \sum_{n=1}^{\infty}\left(\frac{1}{n}+(-1)^{n}\frac{1}{n2^{n+1}}\right)^{2}\\
&=& \frac{\pi^{2}}{6}+ \frac{1}{4}{\rm Li\,}_{2}\left(\frac{1}{4}\right) +{\rm Li\,}_{2}\left(-\frac{1}{2}\right)\\
&=& \frac{\pi^{2}}{6}+\frac{1}{2}\left [ {\rm Li\,}_{2}\left(\frac{1}{2}\right)+3{\rm Li\,}_{2}\left(\frac{-1}{2}\right)\right ],
\eeqq
where we have used the fact that ${\rm Li\,}_{2}(z^2)=2({\rm Li\,}_{2}(z)+{\rm Li\,}_{2}(-z))$.
From the proof of Theorem \ref{13-th 12b}, we conclude  that
$$\sum_{n=1}^{\infty}|\gamma_{n}(f_{1})|^{2}<\frac{\pi^{2}}{6},
$$
because
$$ {\rm Li\,}_{2}\left(\frac{1}{2}\right)+ 3 {\rm Li\,}_{2}\left(\frac{-1}{2}\right) <0.
$$
As a direct approach, it is easy to see that
$$ {\rm Li\,}_{2}(z)+3{\rm Li\,}_{2}(-z)  = \sum_{n=1}^{\infty} \frac{1}{n^2}(1+3(-1)^n)z^n
= \sum_{k=1}^{\infty} \frac{z^{2k}}{k^2}-2\sum_{k=1}^{\infty} \frac{z^{2k-1}}{(2k-1)^2}
$$
and thus, we obtain that
$$ \sum_{n=1}^{\infty}|\gamma_{n}(f_1)|^{2}
=\frac{\pi^{2}}{6}+\frac{1}{2}\sum_{k=1}^{\infty}\frac{1}{4^k} \left(\frac{1}{k^2}-\frac{1}{(k-1/2)^2} \right)
=\frac{\pi^{2}}{6}-\sum_{k=1}^{\infty}\frac{1}{4^k} \left(\frac{4k-1}{k^2(2k-1)^2} \right)
$$
and thus,
$$ \sum_{n=1}^{\infty}|\gamma_{n}(f_1)|^{2} <\frac{\pi^{2}}{6}.
$$
On the other hand, it is a simple exercise to verify that  $f_1 \notin \mathcal{S}^{\star}$. The graph of this function is shown in
Figure \ref{logCoeffig1}(d).
\er

Let ${\mathcal G}(\alpha )$ denote the class of locally univalent normalized analytic
functions $f$ in the unit disk $|z| < 1$ satisfying the condition
$${\rm Re} \left (1+\frac{zf''(z)}{f'(z)}\right )<1+\frac{\alpha}{2}
\quad \mbox{for $|z|<1$,}
$$
and for some $0<\alpha \leq 1$. Set ${\mathcal G}(1)=:{\mathcal G}$. It is known  (see \cite[Equation (16)]{PoRa1}) that
$\mathcal{G}\subset {\mathcal S}^{\star}$ and thus, functions in ${\mathcal G}(\alpha )$ are starlike. This class has been studied
extensively in the recent past, see for instance \cite{OPW-2013} and the references therein. We now consider the estimate of the
type \eqref{eq36} for the subclass ${\mathcal G}(\alpha )$.

\bthm\label{13-th 23b}
Let $0<\alpha \leq 1$ and ${\mathcal G}(\alpha )$ be defined as above. Then
the logarithmic coefficients $\gamma_{n}$ of $f\in {\mathcal G}(\alpha )$ satisfy the inequalities
\be\label{eq 85}
\sum_{n=1}^{\infty}n^{2}|\gamma_{n}|^{2}\leq\frac{\alpha }{4(\alpha +2)}
\ee
and
\be\label{eq 85d}
\sum_{n=1}^{\infty}|\gamma_{n}|^{2}\leq \frac{\alpha ^2}{4}\,{\rm Li\,}_{2}\left(\frac{1}{(1+\alpha )^2} \right ).
\ee
%Equality is attained for the function $f$ such that $f'(z)=(1-z)^{\alpha }$.
Also we have
\be\label{eq 85a}
|\gamma_{n}|\leq \frac{\alpha }{2(\alpha +1)n} ~\mbox{ for $n\ge 1$}.
\ee
\ethm
\bpf
If $f\in {\mathcal G}(\alpha)$, then we have (see eg.  \cite[Theorem 1]{JO-95} and \cite{PoRa1})
\be\label{eq 86}
\frac{zf'(z)}{f(z)} -1 \prec  \frac{(1+\alpha )(1-z)}{1+\alpha -z}-1 = -\alpha \left (\frac{z/(1+\alpha)}{1-(z/(1+\alpha))}\right ), \quad \mbox{$z\in {\mathbb D}$},
\ee
which, in terms of the logarithmic coefficients $\gamma_{n}$ of $f$ defined by \eqref{eq:log}, is equivalent to
\be\label{eq 87}
\sum_{n=1}^{\infty}(-2n \gamma_{n})z^{n}\prec  \alpha\sum_{n=1}^{\infty}\frac{z^n}{(1+\alpha )^{n}}.
\ee
Again, by Rogosinski's result, we obtain that
$$ \sum_{n=1}^{\infty}4n^{2} |\gamma_{n}|^{2}\leq \alpha ^2\sum_{n=1}^{\infty}\frac{1}{(1+\alpha )^{2n}}=\frac{\alpha }{\alpha +2}
$$
which is \eqref{eq 85}.

%The function $f$ such that $f'(z)=(1-z)^{\alpha }$ shows that the inequalities \eqref{eq 85} and \eqref{eq 85d} are best possible.

Now, since the sequence
$A_{n}=\frac{1}{(1+\alpha)^{n}}$ is convex decreasing,  we obtain from \eqref{eq 87} and \cite[Theorem VII, p.64]{Rogo43} that
$$|-2n \gamma_{n}|\leq \alpha A_{1}=\frac{\alpha }{1+\alpha},
$$
which implies the desired inequality \eqref{eq 85a}. As an alternate approach to prove this inequality,  we may rewrite \eqref{eq 86} as
$$\sum_{n=1}^{\infty} (2n \gamma_{n})z^{n} =z\left (\log\left (\frac{f(z)}{z}\right )\right )' \prec \phi (z)=-\alpha \left (\frac{z/(1+\alpha)}{1-(z/(1+\alpha))}\right )
$$
and, since $\phi (z)$ is convex in $\ID$ with $\phi '(0)=-\alpha /(1+\alpha)$, it follows from Rogosinski's result
(see also \cite[Theorem 6.4(i), p.~195]{Duren}) that $|2n \gamma_{n}|\leq \alpha /(1+\alpha)$. Again, this proves the inequality \eqref{eq 85a}.

Finally, we  prove the inequality \eqref{eq 85d}. From the formula \eqref{eq 87} and the
result of Rogosinski (see also \cite[Theorem 2.2]{P} and \cite[Theorem 6.2]{Duren}), it follows that for $k\in \IN$ the inequalities
$$\sum_{n=1}^kn^2\left|\gamma_n\right|^2 \leq \frac{\alpha ^2}{4}\sum_{n=1}^k \frac{1}{(1+\alpha)^{2n}}
$$
are valid. Clearly, this implies the inequality \eqref{eq 85} as well. On the other hand, consider these inequalities for $k=1, \ldots ,N$,
and multiply the $k$-th inequality by the factor $\frac{1}{k^2}-\frac{1}{(k+1)^2},$ if $k=1, \ldots ,N-1$ and by
$\frac{1}{N^2}$ for $k=N$. Then the summation of the multiplied inequalities yields
\beqq
\sum_{k=1}^N\left|\gamma_k\right|^2 &\leq & \frac{\alpha ^2}{4} \sum_{k=1}^N\frac{1}{k^2 (1+\alpha)^{2k}}\\
& \leq & \frac{\alpha ^2}{4}\sum_{k=1}^{\infty}\frac{1}{k^2 (1+\alpha)^{2k}}\\
&=&\frac{\alpha ^2}{4} \,{\rm Li\,}_{2}\left(\frac{1}{(1+\alpha )^2} \right ) ~\mbox{ for $N=1,2,\ldots ,$}
\eeqq
which proves the desired assertion \eqref{eq 85d} if we allow $N\rightarrow \infty$.
%The choice $f$ such that $f'(z)=(1-z)^{\alpha }$ shows that the inequality is best possible.
\epf

\bcor\label{13-th 23}
The logarithmic coefficients $\gamma_{n}$ of $f\in \mathcal{G}:=\mathcal{G}(1)$ satisfy the inequalities
$$\sum_{n=1}^{\infty}n^{2}|\gamma_{n}|^{2}\leq\frac{1}{12} ~\mbox{ and }~
 \sum_{n=1}^{\infty}|\gamma_{n}|^{2}\leq \frac{1}{4}\,{\rm Li\,}_{2}\left(\frac{1}{4} \right ).
$$
The results are  the best possible as the function $f_0(z)=z-\frac{1}{2}z^{2}$ shows. Also we have
$|\gamma_{n}|\leq 1/(4n)$ for $n\ge 1$.
\ecor

\br
For the function $f_0(z)=z-\frac{1}{2}z^{2}$, we have that $\gamma_{n}(f_0)= -\frac{1}{n2^{n+1}}$ for $n=1,2,\ldots$
and thus, it is reasonable to expect that the inequality $|\gamma_{n}|\leq \frac{1}{n2^{n+1}}$ is valid for
the logarithmic coefficients $\gamma_{n}$ of each $f\in \mathcal{G}.$
But that is not the case as the function $f_{n}$ defined by $f_{n}'(z)=(1-z^{n})^{\frac{1}{n}}$ shows. Indeed for this function we have
$$1+\frac{zf_n''(z)}{f_n'(z)} =\frac{1-2z^n}{1-z^n}
$$
showing that $f_n\in \mathcal{G}$. Moreover,
$$\log \frac{f_{n}(z)}{z}=-\frac{1}{n(n+1)}z^{n}+ \cdots,
$$
which implies that $|\gamma_{n}(f_n)|= \frac{1}{2n(n+1)}$ for $n=1,2,\ldots$, and observe that
$\frac{1}{2n(n+1)}>\frac{1}{n2^{n+1}}$ for $n=2,3, \ldots$. Thus, we conjecture that the logarithmic coefficients
$\gamma_{n}$ of each $f\in \mathcal{G}$ satisfy the inequality $|\gamma_{n}|\leq\frac{1}{2n(n+1)}$ for $ n=1,2,\ldots $.
Clearly, Corollary \ref{13-th 23} shows that the conjecture is true for $n=1$.
\er

\br\label{rem3}
Let $f\in \mathcal{C} (\alpha )$, where $0\leq \alpha <1$. Then  we have \cite{WF-80}
\be \label{eq 87f}
\frac{zf'(z)}{f(z)}-1 \prec G_\alpha (z)-1 =\sum_{n=1}^{\infty}\delta_n z^n,
\ee
where $\delta_n$ is real for each $n$,
$$  G_\alpha (z)= \left \{ \begin{array}{ll}
\ds\frac{(2\alpha -1)z}{(1-z)[(1-z)^{1-2\alpha} -1]} &\mbox{ if $\alpha \neq 1/2$,}\\[4mm]
\ds \frac{-z}{(1-z)\log (1-z)} &\mbox{ if $\alpha = 1/2$,}
\end{array}
\right .
$$
and
$$\beta (\alpha)=G_{\alpha}(-1)=\inf_{|z|<1} G_{\alpha}(z)=
\left \{ \begin{array}{ll}
\ds\frac{1-2\alpha}{2[2^{1-2\alpha} -1]} &\mbox{ if $0\leq \alpha \neq 1/2<1$,}\\[4mm]
\ds \frac{1}{2\log 2} &\mbox{ if $\alpha = 1/2$}
\end{array}
\right .
$$
so that $f\in {\mathcal S}^{\star}(\beta (\alpha))$.
Also, we have \cite{Silvia85}
$$\frac{f(z)}{z}\prec  \frac{K_\alpha (z)}{z}= \left \{ \begin{array}{ll}
\ds\frac{(1-z)^{2\alpha -1} -1}{(1-2\alpha)z} &\mbox{ if $0\leq \alpha \neq 1/2<1$,}\\[4mm]
\ds -\frac{\log (1-z)}{z} &\mbox{ if $\alpha = 1/2$,}
\end{array}
\right .
$$
and $K_\alpha (z)/z$ is univalent and convex (not normalized in the usual sense) in $\ID$.
%Moreover,
%$$\frac{zK_\alpha '(z)}{K_\alpha (z)} =G_\alpha (z).
%$$

Now, the subordination relation \eqref{eq 87f}, in terms of the logarithmic coefficients $\gamma_{n}$ of $f$ defined by \eqref{eq:log}, is equivalent to
$$2\sum_{n=1}^{\infty}n \gamma_{n}z^{n}\prec   G_\alpha (z)-1 = \sum_{n=1}^{\infty}\delta_n z^n, \quad z\in\ID,
$$
and thus,
\be\label{eq 87g}
\sum_{n=1}^kn^2\left|\gamma_n\right|^2 \leq \frac{1}{4}\sum_{n=1}^k \delta _n^2 \quad \mbox{ for each $k\in \IN$}.
\ee
Since $f$ is starlike of order $\beta (\alpha)$, it follows that
$$ \frac{zK_\alpha '(z)}{K_\alpha (z)}-1= G_\alpha (z)-1 \prec 2(1-\beta (\alpha))\frac{z}{1-z}
$$
and therefore, $|\delta _n|\leq 2(1-\beta (\alpha))$ for each $n\geq 1$. Again, the relation \eqref{eq 87g} by the previous approach gives
$$\sum_{k=1}^N\left|\gamma_k\right|^2 \leq  \frac{1}{4} \sum_{k=1}^N\frac{\delta _k^2}{k^2}
\leq  (1-\beta (\alpha))^2 \sum_{k=1}^N\frac{1}{k^2}
$$
for $N=1,2,\ldots ,$ and hence, we have
%$$\sum_{n=1}^{\infty}|\gamma_{n}|^{2}\leq (1-\beta (\alpha))^2 \frac{\pi^2}{6}
%$$
%which is not sharp unless $\alpha =0$. However, using \eqref{eq 87g}, one can obtain the following sharp inequality
%$$\sum_{n=1}^{\infty}|\gamma_{n}|^{2}\leq \frac{1}{4}\sum_{n=1}^{\infty} \frac{\delta _n^2}{n^2}
%$$
$$\sum_{n=1}^{\infty}|\gamma_{n}|^{2}\leq \frac{1}{4}\sum_{n=1}^{\infty} \frac{\delta _n^2}{n^2}\leq
 (1-\beta (\alpha))^2 \sum_{n=1}^{\infty} \frac{1}{n^2} =
 (1-\beta (\alpha))^2 \frac{\pi^2}{6}
$$
and equality holds in the first inequality for $K_\alpha (z)$. In particular, if $f$ is convex then $\beta (0)=1/2$ and hence, the last inequality reduces to
$$\sum_{n=1}^{\infty}|\gamma_{n}|^{2}\leq \frac{\pi^2}{24}
$$
which is sharp as the convex function $z/(1-z)$ shows.
\er

 \section{Proof of Conjecture \ref{conj1} for $n=2,3, 4$\label{sec3}}
\bthm\label{13-th 24}
Let $f\in \mathcal{U}(\lambda)$ for $0<\lambda\leq 1$ and let $f(z)=z+a_{2}z^{2}+a_{3}z^{3}+\cdots$.  Then
\be\label{eq 88}
|a_{n}|\leq\frac{1-\lambda^{n}}{1-\lambda} ~\mbox{ for $0<\lambda <1$ and $n=2,3, 4$,}
\ee
and $|a_{n}|\leq n $ for $\lambda=1 $ and $n\geq 2$. The results are the best possible.
\ethm\bpf
The case $\lambda=1 $ is well-known because $\mathcal{U}=\mathcal{U}(1)\subset \mathcal{S}$ and hence, by the de Branges theorem, we have
$|a_{n}|\leq n$ for  $f\in \mathcal{U}$ and $n\geq 2$. Here is an alternate proof without using the de Branges theorem.
From the subordination result \eqref{OPW7-eq4c} with $\lambda=1 $, one has
$$\frac{f(z)}{z}\prec \frac{1}{(1-z)^{2}}=\sum_{n=1}^{\infty} nz^{n-1}
$$
and thus, by  Rogosinski's theorem \cite[Theorem 6.4(ii), p.~195]{Duren}, it follows that  $|a_{n}|\leq n $ for  $n\geq 2$.

So, we may consider $f\in \mathcal{U}(\lambda)$ with $0<\lambda <1$. The result for $n=2$, namely, $|a_{2}|\leq 1+\lambda$ is
proved in \cite{OPW,VY2013} and thus, it suffices to prove \eqref{eq 88} for $n=3, 4$ although our proof below is elegant
and simple for the case $n=2$ as well. To do this, we begin to recall from \eqref{OPW7-eq4c} that
$$\frac{f(z)}{z}\prec \frac{1}{(1-z)(1-\lambda z)}=1+\sum_{n=1}^{\infty} \frac{1-\lambda^{n+1}}{1-\lambda}z^{n}
$$
and thus
$$\frac{f(z)}{z} = \frac{1}{(1-z\omega (z))(1-\lambda z \omega (z))},
$$
where $\omega$ is analytic in $\ID$ and $|\omega(z)|\leq 1$ for $z\in \ID$. In terms of series formulation, we have
$$%\be\label{eq 89}
\sum_{n=1}^{\infty}a_{n+1}z^{n}=\sum_{n=1}^{\infty} \frac{1-\lambda^{n+1}}{1-\lambda}\omega ^n(z)z^n.
$$
We now set $\omega(z)=c_{1}+c_{2}z +\cdots $ and rewrite the last relation as
\be\label{eq 90}
\sum_{n=1}^{\infty}(1-\lambda)a_{n+1}z^{n}=\sum_{n=1}^{\infty} (1-\lambda^{n+1})(c_{1}+c_{2}z +\cdots )^{n}z^{n}.
\ee
By comparing the coefficients of $z^{n}$ for $n=1,2,3$ on both sides of \eqref{eq 90}, we obtain
\be\label{eq 91}
 \left\{\begin{array}{ll}
(1-\lambda)a_{2} & = (1-\lambda^{2})c_{1}\\
(1-\lambda)a_{3}& = (1-\lambda^{2})c_{2}+(1-\lambda^{3})c_{1}^{2}\\
(1-\lambda)a_{4} % & = (1-\lambda^{2})c_{3}+2(1-\lambda^{3})c_{1}c_{2}+(1-\lambda^{4})c_{1}^{3}\\
& = (1-\lambda^{2})\left (c_{3}+\mu c_{1}c_{2}+\nu c_{1}^{3}\right ),
\end{array}
\right.
\ee
where
$$\mu = 2\frac{1-\lambda^3}{1-\lambda^2}  ~\mbox{ and }~ \nu = \frac{1-\lambda^4}{1-\lambda^2}.
$$
It is well-known that $|c_1|\leq 1$ and $|c_{2}|\leq 1-|c_{1}|^{2}$.
From the first relation in \eqref{eq 91} and the fact that $|c_1|\leq 1$, we obtain
$$(1-\lambda)|a_{2}|=(1-\lambda^{2})|c_{1}|\leq 1-\lambda^{2},
$$
which gives a new proof for the inequality $|a_{2}|\leq 1+\lambda$.

Next we present a proof of \eqref{eq 88} for $n=3$. Using the second relation in \eqref{eq 91}, $|c_1|\leq 1$
and the inequality $|c_{2}|\leq 1-|c_{1}|^{2}$, we get
\beqq
(1-\lambda)|a_{3}|&\leq&(1-\lambda^{2})|c_{2}|+(1-\lambda^{3})|c_{1}|^{2}\\
&\leq& (1-\lambda^{2})(1-|c_{1}|^{2})+(1-\lambda^{3})|c_{1}|^{2}\\
%&\leq& (1-\lambda^{3})|c_{1}|^{2}\\
&=& 1-\lambda^{2}+(\lambda^{2}-\lambda^{3})|c_{1}|^{2}\\
%$&\leq&1-\lambda^{2}+\lambda^{2}-\lambda^{3}\\
&\leq &1-\lambda^{3},
\eeqq
which implies $|a_3|\leq 1+\lambda +\lambda ^2$.

Finally, we present a proof of \eqref{eq 88} for $n=4$.  To do this, we recall the sharp upper bounds for  the functionals
$ \left|c_3+\mu c_1c_2+\nu c_1^3\right|$ when  $\mu$ and $\nu $  are real. In \cite{PS}, Prokhorov and Szynal  proved among other results that
$$ \left|c_3+\mu c_1c_2+\nu c_1^3\right|\leq  |\nu|
$$
if $ 2\leq|\mu|\leq 4$ and $\nu \geq (1/12)(\mu^2+8)$. From the third relation in \eqref{eq 91}, this condition is fulfilled and thus,
we find that
$$%\be\label{eq 92a}
(1-\lambda)|a_{4}| =(1-\lambda^{2})\left |c_{3}+\mu c_{1}c_{2}+\nu c_{1}^{3}\right |\leq (1-\lambda^{2})\left (\frac{1-\lambda^4}{1-\lambda^2}\right ) =1-\lambda^4
$$
which proves the desired inequality $|a_4|\leq 1+\lambda +\lambda ^2+\lambda ^3$.
\epf

%We conclude the paper with the following conjecture.
%
%\bcon
%Functions in $\mathcal{U}$ are convex in some direction.
%\econ

\subsection*{Acknowledgements}
%he authors thank the referee for his/her careful reading and many useful comments.
The work of the first author was supported by MNZZS Grant, No. ON174017, Serbia.
The second author is on leave from the IIT Madras.

\end{document}